\input amstex
\magnification=\magstep1 
\baselineskip=13pt
\documentstyle{amsppt}
\vsize=8.7truein \CenteredTagsOnSplits \NoRunningHeads
\def\conv{\operatorname{conv}}

\def\spa{\operatorname{span}}

\def\sym{\operatorname{Sym}}

 \topmatter

\title  Thrifty approximations of convex bodies by polytopes \endtitle
\author Alexander Barvinok  \endauthor
\address Department of Mathematics, University of Michigan, Ann Arbor,
MI 48109-1043, USA \endaddress
\email barvinok$\@$umich.edu \endemail
\date July 2012 \enddate
\thanks  This research was partially supported by NSF Grant DMS 0856640.
\endthanks 
\abstract Given a convex body $C \subset {\Bbb R}^d$ containing the origin in its interior and a real number $\tau > 1$ 
we seek to construct a polytope $P \subset C$ 
with as few vertices as possible such that $C \subset \tau P$. Our construction is nearly optimal for a wide range of 
$d$ and $\tau$. In particular, we prove that 
if $C=-C$ then for any $1>\epsilon>0$ and $\tau=1+\epsilon$ one can choose $P$ having roughly 
$\epsilon^{-d/2}$ vertices and for $\tau=\sqrt{\epsilon d}$ one can choose $P$ having roughly 
$d^{1/\epsilon}$ vertices. Similarly, we prove that if $C \subset {\Bbb R}^d$ is a convex body such that 
$-C \subset \mu C$ for some $\mu \geq 1$ then one can choose 
$P$ having roughly $\bigl((\mu+1)/(\tau-1)\bigr)^{d/2}$ vertices provided  $(\tau-1)/(\mu+1) \ll 1$.
\endabstract
\endtopmatter

\document

\head 1. Introduction and main results \endhead

We discuss how well convex bodies (compact convex sets with non-empty interior) can be approximated by 
polytopes (convex hulls of finite sets of points). There is, of course, a vast literature on the topic, as there are many
different notions of approximation, see surveys \cite{G93a} and \cite{Br07}.
Our setup is as follows. Let 
$C \subset {\Bbb R}^d$ be a convex body containing the origin in its interior. We seek to construct a polytope 
$P \subset {\Bbb R}^d$ with as few vertices as possible, so that 
$$P \subset C \subset \tau P$$
for some given $\tau >1$.

Our first main result concerns {\it symmetric} convex bodies $C$ for which $C=-C$ and $\tau$ measures
the Banach-Mazur distance between $P$ and $C$.

\proclaim{(1.1) Theorem} Let $k$ and $d$ be positive integers and let $\tau > 1$ be a real number such that 
$$\left(\tau-\sqrt{\tau^2 -1}\right)^k +\left(\tau+\sqrt{\tau^2-1}\right)^k \ \geq \ 6 {d+k \choose k}^{1/2}.$$
Then for any symmetric convex body $C \subset {\Bbb R}^d$ there is a symmetric polytope $P \subset {\Bbb R}^d$ with 
at most 
$$8 {d+k \choose k}$$ 
vertices such that 
$$P \ \subset \ C \ \subset \ \tau P.$$
\endproclaim

In fact, we can replace $\displaystyle {d+k \choose k}$ throughout the statement of Theorem 1.1 by a slightly smaller number
$$D(d, k)=\sum_{m=0}^{\lfloor k/2 \rfloor} {d+k-1-2m \choose k-2m}. \tag1.1.1$$
For example, taking $d=20$ and $k=3$ we conclude that any $20$-dimensional symmetric convex body 
can be approximated within a factor of $\tau=3.18$ by a symmetric polytope with at most $12,480$ vertices.

Taking $\tau$ in Theorem 1.1 arbitrarily close to 1, we obtain the following corollary.

\proclaim{(1.2) Corollary} For any 
$$ \gamma \ > \ {e \over 4 \sqrt{2}} \approx 0.48$$ 
there exists $\epsilon_0=\epsilon_0(\gamma) >0$ such that for any $0 < \epsilon < \epsilon_0$ and for any 
symmetric convex body $C \subset {\Bbb R}^d$, there is a symmetric polytope $P \subset {\Bbb R}^d$
with at most 
$$\left({ \gamma  \over \sqrt{\epsilon}} \ln {1 \over \epsilon} \right)^d$$
vertices 
such that 
$$P \ \subset \ C \ \subset \ (1+\epsilon)P.$$ 
\endproclaim

The well-known volumetric argument (see, for example, Lemma 4.10 of \cite{Pi89}) produces polytopes 
with roughly $(3/\epsilon)^d$ vertices which approximate a given symmetric $d$-dimensional convex body within a factor of 
$1+\epsilon$. Hence for small $\epsilon > 0$ the estimate of Corollary 1.2 gives us roughly the square root of the number 
of vertices required by the volumetric bound. It follows from results of Dudley \cite{Du74} and also from results of 
Bronshtein and Ivanov \cite{BI75} that in any dimension $d$
one can construct a polytope $P$ with not more than $\gamma(d) \epsilon^{-(d-1)/2}$ vertices approximating a given symmetric convex body $C \subset {\Bbb R}^d$ within a factor of $1+\epsilon$, with $\gamma(d)$ of the order of $d^{d/2}$.
If the boundary of $C$ is ${\Cal C}^2$-smooth then for all sufficiently small $0<\epsilon < \epsilon_0(C)$ one can obtain an 
approximating symmetric polytope with at most $(\gamma/ \epsilon)^{(d-1)/2}$ vertices for some absolute constant $\gamma>0$, and the dependence on $\epsilon$ cannot be made better \cite{G93b}, \cite{B\"o00} (note that the upper 
bound for $\epsilon$ depends on the convex body $C$).
The estimate of Corollary 1.2 is the first bound improving the volumetric bound uniformly over all symmetric convex bodies $C$ of all dimensions $d$.

Next, we consider approximations for which we want to keep the number of vertices of the polytope polynomial in the dimension
of the ambient space.
\proclaim{(1.3) Corollary} For any 
$$\gamma \ > \ {\sqrt{e} \over 2} \approx 0.82$$ there is 
a positive integer $k_0=k_0(\gamma)$ such that for any $k > k_0$ and for any symmetric convex body $C \subset {\Bbb R}^d$ of a sufficiently large dimension $d > d_0(k)$
there is a symmetric polytope $P \subset {\Bbb R}^d$
with at most 
$$8 {d+k \choose k} $$ 
vertices such that 
$$P \ \subset \ C \ \subset \ \gamma \sqrt{d \over k} P.$$
\endproclaim

A simple computation shows that if $C$ is the $d$-dimensional Euclidean ball and $P$ has at most $d^k$ vertices for some fixed 
$k$, then $P$ cannot approximate $C$ better than within a factor of $\displaystyle \tau= \gamma \sqrt{d \over k \ln d}$ as 
$d$ grows, where $\gamma>0$ is an absolute constant.

Finally, we consider approximations of not necessarily symmetric convex bodies. We prove the following main result,
generalizing Theorem 1.1. The quality of approximation depends on the {\it symmetry coefficient} of the convex body $C$,
that is on the smallest $\mu \geq 1$ such that $-C \subset \mu C$ (recall that the convex bodies we consider contain the 
origin in their interior). 

\proclaim{(1.4) Theorem} Let $d$ and $k$ be positive integers. For $\tau, \mu \geq 1$ let us define 
$$\lambda=\lambda(\tau, \mu)={2 \over \mu+1} \tau + {\mu-1 \over \mu+1} \ \geq\ 1.$$
If 
$$\left(\lambda -\sqrt{\lambda^2-1}\right)^k +\left(\lambda+\sqrt{\lambda^2-1}\right)^k \ \geq \ 6{d+k \choose k}^{1/2}$$
then for any convex body $C \subset {\Bbb R}^d$ containing the origin in its interior and such that 
$-C \subset \mu C$ there is a polytope $P \subset {\Bbb R}^d$ with  
at most
$$8{d+k \choose k}$$
vertices such that 
$$P \ \subset \ C \ \subset \ \tau P.$$
\endproclaim
We also obtain the following extension of Corollary 1.2.

\proclaim{(1.5) Corollary} 
\roster
\item
For $\tau, \mu \geq 1$ let us define 
$$\delta=\delta(\tau, \mu)={2(\tau-1) \over \mu+1}.$$ For any 
$$\gamma \ > \ {e \over 4 \sqrt{2}} \approx 0.48$$
there exists $\delta_0=\delta_0(\gamma) >0$ such that as long as $\delta(\tau, \mu)< \delta_0$, for any convex body 
$C \subset {\Bbb R}^d$ such that $-C \subset \mu C$ there exists a polytope with at most 
$$\left({\gamma \over \sqrt{\delta}} \ln {1 \over \delta}\right)^d$$
vertices such that 
$$P \ \subset \ C \ \subset \ \tau P.$$
\item For any 
$$\gamma \ > \ {e \over 8} \approx 0.34$$
there exists $\epsilon_0=\epsilon(\gamma)>0$ such that for any $0 < \epsilon < \epsilon_0$ and for any convex body 
$C \subset {\Bbb R}^d$ such that $-C \subset \mu C$ for some $\mu \geq 1$ there exists a polytope 
$P \subset {\Bbb R}^d$ with at most 
$$\left(\gamma \sqrt{\mu+1 \over \epsilon} \ln {1 \over \epsilon} \right)^d$$
vertices such that 
$$P \ \subset \ C \ \subset \ (1+\epsilon) P.$$
\endroster
\endproclaim
As a function of the symmetry coefficient $\mu$, the number of vertices of $P$ grows roughly as $\mu^{d/2}$ as long as
the ratio $\tau/\mu$ is small enough. It follows from results of Gruber \cite{G93b} that if the boundary of $C$ is ${\Cal C}^2$-smooth
then for all sufficiently small $0<\epsilon < \epsilon_0(C)$ one can construct a polytope $P$ with not more than 
$\mu^{d/2} (\gamma/\epsilon)^{(d-1)/2}$ vertices for some absolute constant $\gamma$ which approximates $C$ within a factor of $1+\epsilon$. The estimates of Corollary 1.5 are uniform over all convex bodies $C$ of all dimensions $d$.

The plan of the paper is as follows. In Section 2, we collect some facts needed for the proofs of Theorems 1.1 and 1.4.
Namely, we review the classical result on the John decomposition of the identity operator and the minimum volume ellipsoid 
of a convex body, a recent result of Batson, Spielman and Srivastava \cite{B+08} which allows one to obtain certain ``sparsification" of the John decomposition, the standard construction of tensor product from multilinear 
algebra which allows us to translate polynomial relations among vectors into linear identities among tensors and
the classical construction of the Chebyshev polynomials which solve a relevant extremal problem. As it turns out, the vertices of the approximating polytopes $P$ are picked up by certain algebraic conditions.

We complete the proofs in Section 3. 

\head 2. Preliminaries \endhead

\subhead (2.1) Chebyshev polynomials \endsubhead
For a positive integer $k$ let $T_k(t)$ be the Chebyshev polynomial of degree $k$, see, for example,
Section 2.1 of \cite{BE95}. Thus for real $t$ the polynomial
$T_k(t)$ can be defined by 
$$\split &T_k(t)= \cos\left(k \arccos t\right) \quad \text{provided} \quad -1 \leq t \leq 1 \quad \text{and} \\
&T_k(t)={1 \over 2} \left(t -\sqrt{t^2-1}\right)^k +{1 \over 2} \left(t + \sqrt{t^2-1}\right)^k \quad \text{provided}
\quad |t| \geq 1. \endsplit$$
In particular, 
$$\left|T_k(t)\right| \ \leq 1 \quad \text{provided} \quad |t| \leq 1. \tag2.1.1$$
Writing $T_k(t)$ in the standard monomial basis, we obtain
$$T_k(t) ={k \over 2} \sum_{m=0}^{\lfloor k/2 \rfloor} (-1)^m {(k-m-1)! \over m! (k-2m)!} (2t)^{k-2m}.$$
In particular,
$$T_1(t)=t, \ T_2(t)=2t^2-1,\ T_3(t)=4t^3-3t,\ T_4(t)=8t^4-8t^2+1.$$
We note that $T_k(-t)=T_k(t)$ if $k$ is even and $T_k(-t)=-T_k(t)$ if $k$ is odd. We also note that the polynomial
$T_k(t)$ is strictly increasing for $t \geq 1$.

In particular,
$$\left| T_k(t) \right| \ > \ { \left(\tau -\sqrt{\tau^2-1}\right)^k + \left(\tau + \sqrt{\tau^2-1}\right)^k  \over 2} \quad 
\text{provided} \quad |t| > \tau \geq 1. \tag2.1.2$$
The polynomial $T_k(t)$ has the following extremal property relevant to us: for any $t_0 \notin [-1, 1]$ the maximum 
value of $\left|p(t_0)\right|$, where $p$ is a polynomial of $\deg p \leq k$ such that $\left|p(t)\right| \leq 1$ for all $t \in [-1, 1]$, 
is attained for $p=T_k$, see, for example, Section 5.1 of \cite{BE95}.

\subhead (2.2) Tensor power \endsubhead
Let $V$ be Euclidean space with scalar product $\langle \cdot, \cdot \rangle$. For a positive integer $k$ let 
$$V^{\otimes k}=\underbrace{V \otimes \cdots \otimes V}_{\text{$k$ times}}$$
be the $k$-th tensor power of $V$. We consider $V^{\otimes k}$ as Euclidean space endowed with scalar 
product $\big\langle \cdot, \cdot \big\rangle$ such that 
$$\big\langle x_1 \otimes \cdots \otimes x_k,\ y_1 \otimes \cdots \otimes y_k \big\rangle=
\prod_{i=1}^k \langle x_i, y_i \rangle$$
for all $x_1, \ldots, x_k; y_1, \ldots, y_k \in V$.
The space $V^{\otimes 2}$ is naturally identified with the space of all linear operators on $V$. 

The symmetric part $\sym\left(V^{\otimes k}\right)$ of $V^{\otimes k}$ is the subspace spanned by the
tensors
$$x^{\otimes k}=\underbrace{x \otimes \cdots \otimes x}_{\text{$k$ times}}$$
for $x \in V$. The space $\sym\left(V^{\otimes k}\right)$ is naturally identified with the space of all homogeneous 
polynomials $p: V \longrightarrow {\Bbb R}$ of degree $k$. In particular, $\sym\left(V^{\otimes 2}\right)$ can be 
identified with the space of quadratic forms on $V$ and also with the space of all symmetric operators on $V$.
We have 
$$\dim \sym\left(V^{\otimes k}\right)={\dim V +k -1 \choose k}.$$
Let us consider the direct sum 
$$ W={\Bbb R} \oplus V \oplus V^{\otimes 2} \oplus \cdots \oplus V^{\otimes k}$$
as Euclidean space with the standard scalar product, which we also denote by $\big\langle \cdot , \cdot \big\rangle$.
For a real univariate polynomial $a(t)$
and a vector $x \in V$, we denote by $a^{\otimes}(x) \in W$ the vector 
$$a^{\otimes}(x) =\alpha_0 \oplus \alpha_1 x \oplus \alpha_2 x^{\otimes 2} \oplus \ldots \oplus \alpha_k x^{\otimes k},
\quad \text{where} \quad a(t)=\sum_{m=0}^k \alpha_m t^m. \tag2.2.1$$
It is then easy to check that for any $x, y \in V$ and any polynomials $a(t)$ and $b(t)$, we have 
$$\aligned &\big\langle a^{\otimes}(x),\ b^{\otimes}(y)\big\rangle =c\left(\langle x, y \rangle\right), \quad
\text{provided} \\
&a(t)=\sum_{m=0}^k \alpha_m t^m,\ b(t)=\sum_{m=0}^k \beta_m t^m \quad \text{and} \quad c(t)=\sum_{m=0}^k \left(\alpha_m \beta_m\right) t^m. \endaligned
 \tag2.2.2$$

\subhead (2.3) The ellipsoid of the minimum volume \endsubhead
As is known, for any compact set $C \subset {\Bbb R}^d$ 
there is a unique ellipsoid of the minimum volume among all ellipsoids centered 
at the origin and containing $C$. If the minimum volume ellipsoid is the unit ball
$$B=\bigl\{x \in {\Bbb R}^d: \quad \|x\| \leq 1 \bigr\},$$
where $\| \cdot \|$ is the Euclidean norm, the contact points $x_i \in C \cap \partial B$ provide a certain 
decomposition of the identity operator $I$, called the John decomposition (recall that $x \otimes x$ for $x \in {\Bbb R}^d$ 
is interpreted as a $d \times d$ symmetric matrix). We need the following result, see for example, \cite{Ba97}.

\proclaim{(2.3.1) Theorem} Let $C \subset {\Bbb R}^d$ be a compact set which spans ${\Bbb R}^d$
and let $B \subset {\Bbb R}^d$ be the unit ball.
Suppose that $C \subset B$ and that $B$ has the smallest volume among all ellipsoids centered at the origin and 
containing $C$. Then there exist points $x_1, \ldots, x_n \in C \cap \partial B$ and non-negative real 
$\alpha_1, \ldots, \alpha_n$ such that 
$$\sum_{i=1}^n \alpha_i \left(x_i \otimes x_i \right) =I,$$
where $I$ is the identity operator on ${\Bbb R}^d$. Equivalently,
$$\sum_{i=1}^n \alpha_i \langle x_i, y \rangle^2 = \|y\|^2$$
for every $y \in {\Bbb R}^d$.
\endproclaim

\subhead (2.4) Sparsification \endsubhead We need a recent result of Batson, Spielman and Srivastava on 
a certain ``sparsification" of the conclusion of Theorem 2.3.1. Namely, we want to be able to choose the number $n$ 
of points in Theorem 2.3.1 linear in the dimension $d$ at the cost of a controlled corruption of the identity operator $I$.

If $A$ and $B$ are $d \times d$ symmetric matrices we say that $A \preceq B$ if $B-A$ is positive semidefinite.
The following result is from \cite{B+08}.

\proclaim{(2.4.1) Theorem} Let $\gamma >1$ be a number and let $x_1, \ldots, x_n$ be vectors in ${\Bbb R}^d$ such that 
$$\sum_{i=1}^n x_i \otimes x_i =  I,$$
or, equivalently,
$$\sum_{i=1}^n \langle x_i, y\rangle^2 =\|y\|^2$$
for all $y \in {\Bbb R}^d$.
Then there is a subset $J \subset \{1, \ldots, n\}$ with $|J| \leq \gamma d$ and $\beta_j >0$ for $j \in J$ such that 
$$I \ \preceq \ \sum_{j \in J} \beta_j \left(x_j \otimes x_j\right) \ \preceq \ 
\left({\gamma +1 +2 \sqrt{\gamma} \over \gamma+1 -2 \sqrt{\gamma}}\right) I,$$
or, equivalently,
$$\|y\|^2 \ \leq\ \sum_{j \in J} \beta_j \langle x_j , y\rangle^2 \leq 
\left({\gamma +1 +2 \sqrt{\gamma} \over \gamma+1 -2 \sqrt{\gamma}}\right) \|y\|^2$$
for all $y \in {\Bbb R}^d$.
\endproclaim

\head 3. Proofs \endhead

We start with a lemma.

\proclaim{(3.1) Lemma} Let $C \subset {\Bbb R}^d$ be a compact set. 
Then there is a subset 
$X \subset C$ of 
$$|X| \ \leq \ 4d$$
points such that for any linear function $\ell: {\Bbb R}^d \longrightarrow {\Bbb R}$ we have 
$$\max_{x \in X} | \ell(x)| \ \leq \ \max_{x \in C} | \ell(x)| \ \leq \ 3\sqrt{d} \max_{x \in X} |\ell(x)|.$$
\endproclaim
\demo{Proof} Without loss of generality we assume that $C$ spans ${\Bbb R}^d$.
Applying a linear transformation, if necessary, we may assume that $C$ is contained in the unit ball $B$ 
and that $B$ is the minimum volume ellipsoid among all ellipsoids centered at the origin and containing $C$. 
By Theorem 2.3.1 there exist vectors $x_1,\ldots x_n \in C \cap \partial B$ and numbers 
$\alpha_1, \ldots, \alpha_n \geq 0$ such that
$$\sum_{i=1}^n \alpha_i \left(x_i \otimes x_i\right) = I.$$
Applying Theorem 2.4.1 with $\gamma=4$ to vectors $\sqrt{\alpha_i} x_i$ we conclude that for some 
$J \subset \{1, \ldots, n\}$ and $\beta_j >0$ for $j \in J$ we have
$$I \ \preceq \ \sum_{j \in J} \alpha_j \beta_j \left( x_j \otimes x_j \right) \ \preceq \ 9I \tag3.1.1$$
and $|J| \leq 4d$. We let
$$X=\left\{x_j: \ j \in J \right\}.$$
 In particular, $x_j \in C$ and $\|x_j\|=1$ for all $j \in J$.
Comparing the traces of the operators in (3.1.1), we get
$$d \ \leq \ \sum_{j \in J} \alpha_j \beta_j \ \leq \ 9d. \tag3.1.2$$
A linear function $\ell: {\Bbb R}^d \longrightarrow {\Bbb R}$ can be written as 
$\ell(x)=\langle y,  x \rangle$ for some $y \in {\Bbb R}^d$.  It follows by (3.1.1) that 
$$\sum_{j \in J} \left(\alpha_j \beta_j\right)  \langle y, x_j \rangle^2 \ \geq \  \|y\|^2$$
and then by (3.1.2) it follows that 
$$|\langle y, x _j\rangle| \ \geq \ {1 \over 3\sqrt{d}} \|y\| \quad \text{for some} \quad j \in J.$$
Since $C \subset B$, we have 
$$\max_{x \in C} |\langle y, x \rangle| \ \ \leq \ \|y\|$$
and the proof follows.
{\hfill \hfill \hfill} \qed
\enddemo

We now prove Theorem 1.1.

\subhead (3.2) Proof of Theorem 1.1 \endsubhead
Let us denote $V={\Bbb R}^d$ and let us consider the space 
$$W={\Bbb R} \oplus V \oplus V^{\otimes 2} \oplus \cdots \oplus V^{\otimes k},$$
see Section 2.2.
Let us define a continuous map $\phi: V \longrightarrow W$ by 
$$\phi(x)=1 \oplus x \oplus x^{\otimes 2} \oplus \cdots \oplus x^{\otimes k} \quad \text{for} \quad x \in V.$$
We consider the compact set 
$$\widehat{C}=\left\{\phi(x): \quad x \in C \right\}, \quad \widehat{C} \subset W.$$
We note that $\widehat{C}$ lies in the subspace 
$${\Bbb R} \oplus V \oplus \sym\left(V^{\otimes 2}\right) \oplus \cdots \oplus \sym\left(V^{\otimes k}\right).$$
In particular,
$$\dim \spa\left(\widehat{C}\right) \ \leq \ 1+ d + {d+1 \choose 2} + \ldots + {d+k-1 \choose k}={d+k \choose k}.$$
Applying Lemma 3.1 to $\widehat{C}$, we conclude that there is a set $X \subset C$ such that 
$$|X| \ \leq \ 4 {d+k \choose k}$$ 
such that for any linear function ${\Cal L}: W \longrightarrow {\Bbb R}$ we have 
$$\max_{x \in X} \left|{\Cal L}\bigl(\phi(x)\bigr)\right| \ \leq \ 
\max_{x \in C} \left|{\Cal L}\bigl(\phi(x)\bigr)\right| \ \leq \ 3{d+k \choose k}^{1/2}
 \max_{x \in X} \left|{\Cal L}\bigl(\phi(x)\bigr)\right|. \tag3.2.1$$
We define $P$ as the convex hull 
$$P=\conv\bigl(X \cup -X \bigr).$$
Clearly, $P \subset C$ and $P$ has at most $8{d+k \choose k}$ vertices. To conclude that $P$ approximates $C$ with the
desired accuracy, we compare the maxima of linear functions $\ell: {\Bbb R}^d \longrightarrow {\Bbb R}$ on $C$ and on $P$.

Suppose that 
$$\ell(x)=\langle y, x \rangle \quad \text{for some} \quad y \in V.$$
Let us define a linear function ${\Cal L}: W \longrightarrow {\Bbb R}$ by 
$${\Cal L}(w)=\big\langle T_k^{\otimes}(y),\ w \big\rangle \quad \text{for all} \quad w \in W,$$
where $T_k$ is the Chebyshev polynomial of degree $k$, see Section 2.1 and (2.2.1).
Then by (2.2.2), we have 
$${\Cal L}\bigl(\phi(x)\bigr)=T_k\left( \langle y, x \rangle \right).$$
Hence from (3.2.1) we obtain
$$\max_{x \in X} \left| T_k\bigl(\ell(x)\bigr)\right| \ \leq \ 
\max_{x \in C} \left| T_k\bigl(\ell(x)\bigr)\right| \ \leq \ 3{d+k \choose k}^{1/2}
 \max_{x \in X} \left| T_k\bigl(\ell(x)\bigr)\right|. \tag3.2.2$$
Suppose that  $\ell(x) \leq 1$ for all $x \in P$ and hence $|\ell(x)| \leq 1$ for all $x \in X$.
Then by (2.1.1) we have $\left|T_k\bigl(\ell(x)\bigr)\right| \leq 1$ for all $x \in X$. If for some 
$x \in C$ we have $\ell(x) > \tau$ then by (2.1.2) we have 
$$\left| T_k\left(\ell(x)\right)\right| \ > \ {\left(\tau -\sqrt{\tau^2 -1}\right)^k + \left(\tau +\sqrt{\tau^2-1}\right)^k \over 2} 
\ \geq \ 3{d+k \choose k}^{1/2},$$
which contradicts (3.2.2). Therefore
$$\max_{x \in P} \ell(x) \ \leq \ \max_{x \in C}  \ell(x) \ \leq \ \tau \max_{x \in P} \ell(x) \tag3.2.3$$
for every linear function $\ell: {\Bbb R}^d \longrightarrow {\Bbb R}$,
which proves that $C \subset \tau P$.
{\hfill \hfill \hfill} \qed

\remark{(3.3) Remark} One can sharpen the bounds somewhat by noticing that the polynomial $T_k$ is even for 
even $k$ and odd for odd $k$. Consequently, the map $\phi: V \longrightarrow W$ can be replaced by 
$$\phi_e(x)=1 \oplus x^{\otimes 2} \oplus \cdots \oplus x^{\otimes k-2} \oplus x^{\otimes k}$$
for even $k$ and by 
$$\phi_o(x)=x \oplus x^{\otimes 3} \oplus \cdots  \oplus x^{\otimes k-2} \oplus x^{\otimes k}$$
for odd $k$. This allows us to replace ${d+k \choose k}$ by $D(d, k)$ defined by (1.1.1) throughout the statement 
of Theorem 1.1.
\endremark

\subhead (3.4) Proof of Theorem 1.4 \endsubhead 
As in Section 3.2, we construct the space $W$, the map $\phi$, the set $\widehat{C}$ and the subset $X \subset C$ 
so that (3.2.1) holds. We then define $P$ as the convex hull
$$P=\conv\left(X \cup \left(-1/\mu \right) X \right).$$
Clearly, $P \subset C$ and $P$ has at most $8{d+k \choose k}$ vertices. To conclude that $P$ approximates $C$ with the 
desired accuracy, we compare the maxima of linear functions $\ell: {\Bbb R}^d \longrightarrow {\Bbb R}$ on $C$ and on $P$.

Let $T_k$ be the Chebyshev polynomial of degree $k$. We define a polynomial $S_k$ by 
$$S_k(t)=T_k\left({2 \over \mu+1} t + {\mu-1 \over \mu+1}\right).$$
Hence $\deg S_k(t)=k$. Moreover, 
$$\left| S_k(t) \right|\ \leq \ 1 \quad \text{provided} \quad -\mu \ \leq \ t \  \leq \ 1\tag3.4.1$$
and 
$$\left| S_k(t)\right| \ > \ {\left(\lambda - \sqrt{\lambda^2-1}\right)^k + \left(\lambda +\sqrt{\lambda^2-1}\right)^k \over 2}
\quad \text{provided} \quad t > \tau. \tag3.4.2$$
Given a linear function $\ell: {\Bbb R}^d \longrightarrow {\Bbb R}$, 
$$\ell(x)=\langle y, x \rangle \quad \text{for some} \quad y \in V,$$
we define a linear function ${\Cal L}: W \longrightarrow {\Bbb R}$ by 
$${\Cal L}(w)=\big\langle S_k^{\otimes}(y), \ w \big\rangle \quad \text{for all} \quad w \in W.$$
Then
$${\Cal L}\bigl(\phi(x)\bigr)=S_k\left( \langle y, x \rangle \right).$$
Hence from (3.2.1) we obtain
$$\max_{x \in X} \left| S_k\bigl(\ell(x)\bigr)\right| \ \leq \ 
\max_{x \in C} \left| S_k\bigl(\ell(x)\bigr)\right| \ \leq \ 3{d+k \choose k}^{1/2}
 \max_{x \in X} \left| S_k\bigl(\ell(x)\bigr)\right|. \tag3.4.3$$
Suppose that $\ell(x) \leq 1$ for all $x \in P$. Then, necessarily, $1 \geq \ell(x) \geq -\mu$ for all 
$x \in X$ and hence by (3.4.1) we have 
 $\left|S_k\bigl(\ell(x)\bigr)\right| \leq 1$ for all $x \in X$. If for some 
$x \in C$ we have $ \ell(x)> \tau$ then by (3.4.2)
$$\left| S_k\left(\ell(x)\right)\right| \ > \ {\left(\lambda -\sqrt{\lambda^2 -1}\right)^k + \left(\lambda+\sqrt{\lambda^2-1}\right)^k \over 2} 
\ \geq \ 3{d+k \choose k}^{1/2},$$
which contradicts (3.4.3). Hence (3.2.3) holds for every linear function $\ell: {\Bbb R}^d \longrightarrow {\Bbb R}$
and, therefore, $C \subset \tau P$.
{\hfill \hfill \hfill} \qed

\subhead (3.5) Proof of Corollary 1.2 \endsubhead
Let us choose $\tau=1+\epsilon$ in Theorem 1.1. 
We use the standard estimate 
$${d+k \choose k} \ \leq \ \left({d+k \over k}\right)^k \left({d+k \over d}\right)^d 
 \ \leq \ e^d \left(1 + {k \over d} \right)^d. \tag3.5.1$$
Let us choose
$$k =\left\lceil {\beta d \over \sqrt{\epsilon}} \ln {1 \over \epsilon} \right\rceil, \tag3.5.2$$ 
where $\beta >0$ is a constant.
Then
$$ {1 \over k} \ln\left(6 {d+k \choose d}^{1/2}\right) \ \leq \  {\sqrt{\epsilon} \over 4 \beta}\bigl(1+o(1)\bigr), \tag3.5.3$$ 
where ``$o(1)$" stands for a term which converges to $0$ uniformly on $d$ as $\epsilon \longrightarrow 0$. 

On the other hand,
$$ \ln \left( \tau + \sqrt{\tau^2-1} \right)=\sqrt{2 \epsilon}\left(1+o(1)\right),$$
where ``$o(1)$" stands for a term which converges to $0$ as $\epsilon \longrightarrow 0$.
Then, as long as 
$$\beta \ > \ {1 \over 4 \sqrt{2}},$$
the condition of Theorem 1.1 is satisfied for all sufficiently small $0<\epsilon < \epsilon_0(\beta)$. The proof now follows by
(3.5.1).
{\hfill \hfill \hfill} \qed

\subhead (3.6) Proof of Corollary 1.5 \endsubhead To prove Part (1),
we observe that $\lambda=1+\delta$ and hence
$$ \ln \left( \lambda + \sqrt{\lambda^2-1} \right)=\sqrt{2 \delta}\left(1+o(1)\right), \tag3.6.1$$
where ``$o(1)$" stands for a term which converges to $0$ as $\delta \longrightarrow 0$.
In Theorem 1.4, let us choose $k$ defined by (3.5.2) with $\epsilon$ replaced by $\delta$. Comparing (3.5.3)
with $\epsilon$ replaced by $\delta$ and (3.6.1) we conclude the proof as in 
Section 3.5.

To prove Part (2), in Theorem 1.4 we choose $\tau=1+\epsilon$ and $k$ defined by (3.5.2). Then 
$$\ln\left(\lambda+\sqrt{\lambda^2-1}\right)=2\left({\epsilon \over \mu+1}\right)^{1/2}\Bigl(1+o(1)\Bigr), \tag3.6.2$$
where ``$o(1)$'' stands for a term which converges to 0 uniformly on $\mu \geq 1$ as $\epsilon \longrightarrow 0$.
Comparing (3.6.2) and (3.5.3), we conclude that the condition of Theorem 1.4 is satisfied for all sufficiently small 
$0 < \epsilon < \epsilon_0(\beta)$ as long as 
$$\beta \ > \ {\sqrt{\mu+1} \over 8}.$$
The proof now follows by (3.5.1).
{\hfill \hfill \hfill} \qed

\subhead (3.7) Proof of Corollary 1.3 \endsubhead
Let us choose $\tau=\gamma \sqrt{d/k}$ in Theorem 1.1, where $\gamma >0$ is a constant.
Using Stirling's formula, we conclude that for each $k$
$$\lim_{d \longrightarrow \infty} {1 \over \sqrt{d}} 6^{1/k} {d+k \choose d}^{1/2k} =\sqrt{e \over k} \bigl(1+o(1)\bigr),$$
where ``$o(1)$" stands for a term which converges to $0$ as $k$ grows.

On the other hand, for each $k$
$$\lim_{d \longrightarrow \infty} {\tau + \sqrt{\tau^2 -1} \over \sqrt{d}}={2\gamma \over \sqrt{k}}.$$
The proof now follows by Theorem 1.1.
{\hfill \hfill \hfill} \qed

\head Acknowledgment \endhead

The author is grateful to Mark Rudelson and Roman Vershynin for many helpful conversations.


\Refs
\widestnumber\key{AAAA}

\ref\key{Ba97}
\by K. Ball
\paper An elementary introduction to modern convex geometry
\inbook Flavors of Geometry
\bookinfo Math. Sci. Res. Inst. Publ.
\vol 31
\publ Cambridge Univ. Press
\publaddr Cambridge
\yr1997
\pages 1--58 
\endref

\ref\key{B+08}
\by J. Batson, D.A. Spielman and N. Srivastava
\paper Twice-Ramanujan sparsifiers 
\paperinfo preprint {\tt arXiv:0808.0163}
\yr 2008
\endref

\ref\key{B\"o00}
\by K. B\"or\"oczky
\paper Approximation of general smooth convex bodies
\jour Adv. Math.
\vol 153 
\yr 2000
\pages 325--341
\endref

\ref\key{BE95}
\by P. Borwein and T. Erd\'elyi
\book Polynomials and Polynomial Inequalities
\bookinfo Graduate Texts in Mathematics, 161
\publ Springer-Verlag
\publaddr New York
\yr 1995
\endref

\ref\key{Br07}
\by E.M. Bronshtein
\paper Approximation of convex sets by polyhedra (Russian)
\jour Sovrem. Mat. Fundam. Napravl. 
\vol 22 
\yr 2007
\pages 5--37
\transl  translation in J. Math. Sci. (N. Y.) 153 (2008), no. 6, 727--762 
\endref

\ref\key{BI75}
\by E.M. Bronshtein and L.D. Ivanov
\paper The approximation of convex sets by polyhedra (Russian)
\jour Sibirsk. Mat. Zh. 
\vol 16 
\yr 1975 
\pages 1110--1112
\transl translation in Siberian Math. J. 16 (1975), no. 5, 852--853 (1976)
\endref

\ref\key{Du74}
\by R.M. Dudley
\paper Metric entropy of some classes of sets with differentiable boundaries
\jour J. Approximation Theory 
\vol 10 
\yr 1974 
\pages 227--236
\endref

\ref\key{G93a}
\by P.M. Gruber
\paper Aspects of approximation of convex bodies
\inbook Handbook of Convex Geometry, Vol. A
\pages  319--345
\publ North-Holland
\publaddr Amsterdam
\yr 1993
\endref

\ref\key{G93b}
\by P.M. Gruber
\paper Asymptotic estimates for best and stepwise approximation of convex bodies. I
\jour Forum Math. 
\vol 5 
\yr 1993
\pages 281--297
\endref

\ref\key{Pi89}
\by G. Pisier
\book The Volume of Convex Bodies and Banach Space Geometry
\bookinfo Cambridge Tracts in Mathematics, 94
\publ Cambridge University Press
\publaddr Cambridge
\yr 1989
\endref


\endRefs
\enddocument
\end